\documentclass{amsart}
\usepackage{graphicx} 
\usepackage{amsmath}
\usepackage{amssymb}
\usepackage{amsthm}
\usepackage{amstext}
\usepackage{enumerate}

      \newtheorem{thm}{Theorem}[section]
      \newtheorem{lem}[thm]{Lemma}
      \newtheorem{prop}[thm]{Proposition}
      \newtheorem{cor}[thm]{Corollary}
      \newtheorem{problem}[thm]{Problem}
\theoremstyle{remark}
      \newtheorem{rem}[thm]{\bf Remark}
\title[Signature and crossing number of links]{Signature and crossing number of links}
\author{Kai Ishihara}
\address[K. Ishihara]{Department of Mathematics, Hiroshima University, Higashi-Hiroshima, 739-8526, Japan}
\email{xishihar@hiroshima-u.ac.jp}
\author{Kei Okada}
\address[K. Okada]{Department of Mathematics, Saitama University, Saitama 338-8570, Japan}
\author{Koya Shimokawa}
\address[K. Shimokawa]{Department of Mathematics, Ochanomizu University, Tokyo 112-8610, Japan}
\email{shimokawa.koya@ocha.ac.jp}

\begin{document}

\begin{abstract}
    This paper investigates the relationship between the signature and the crossing number of knots and links. We refine existing theorems and provide a comprehensive classification of links with specific properties, particularly those with signatures that deviate by a fixed amount from their crossing numbers. 
    The main results include the identification of all links for which the sum of the signature and crossing number equals 2, which are shown to be closures of positive 3-braids. 
    Additionally, we explore the implications of these findings in the context of band surgeries and their applications to vortex knots and DNA topology.
\end{abstract}

\maketitle

\section{Introduction}

Let $L$ be an oriented link, which is an embedding of a disjoint union of circles $S^1$ into $S^3$. $L$ is typically depicted using a diagram with arrows that are determined by the orientation of $S^1$. We say $L$ is a knot if it consists of only one component. 
The crossing number, $cr(L)$, is defined as the minimal number of crossings among all diagrams of $L$.
The signature of $L$, $\sigma(L)$, is the signature of the matrix $M+M^T$, where $M$ is the Seifert matrix defined by using a Seifert surface of $L$. 
The concept of the signature was initially introduced for knots by H.F.~Trotter \cite{Trotter} and later extended to arbitrary links by K.~Murasugi \cite{M1}.

Little was known about the relationship between the crossing number and the signature of a link.
In \cite{PNAS}, we give a bound of the signature of an oriented link using its crossing number.

\begin{thm}[\cite{PNAS}]\label{thm:T(2,c)}
Let $L$ be a non-trivial oriented knot or non-split oriented link. 
Then,  
$$-cr(L)<\sigma(L)<cr(L).$$
Moreover, $\sigma(L)=1-cr(L)$ if and only if $L=T(2,c)$, where $c=cr(L)$. 
\end{thm}

\begin{figure}[h]
    \centering
    \includegraphics[scale=.3]{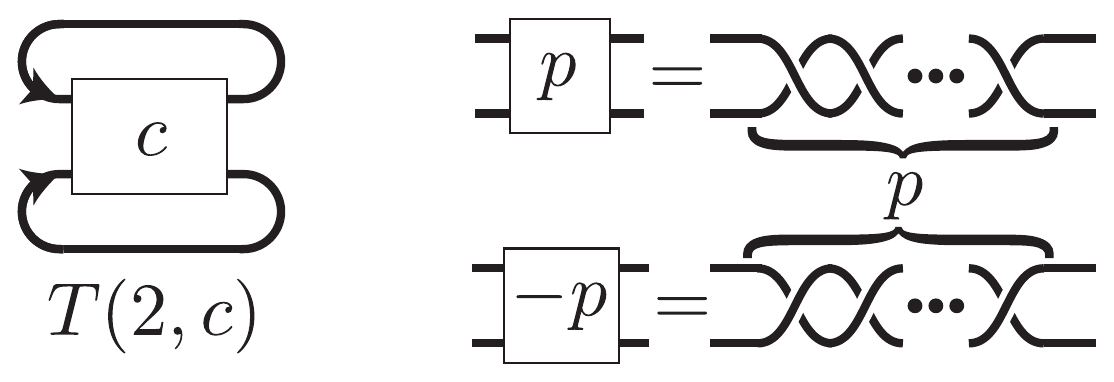}
    \caption{Left: The $(2,c)$-torus link $T(2,c)$ with parallel orientation.
    Right: Each box with a positive integer $p$ (respectively a negative integer $-p$) represents right-handed (resp. left-handed) $p$ twists.}
    \label{fig:2ctorus}
\end{figure}

We use $T(p,q)$ to mean the $(p,q)$-torus link with parallel orientation, see Fig. \ref{fig:2ctorus} for $T(2,c)$. 

We will refine Theorem \ref{thm:T(2,c)} in Section \ref{sec:pre}.
Based on Theorem \ref{thm:T(2,c)}, we can explore Problem \ref{problem}, which is the realization problem of a link that takes specific values for the signature and the crossing number. 
This type of problem is called the geography problem and has been studied by K.~Taniyama for many pairs of knot invariants \cite{Taniyama}.
\begin{problem}\label{problem}
    For any pair of integers $(c,d)$ satisfying $-c<d<c$, is there an oriented link $L$ 
    with $(cr(L),\sigma(L))=(c,d)$?
\end{problem}

It is well known that the mirror image $L!$ of $L$ satisfies $cr(L!)=cr(L)$ and $\sigma(L!)=-\sigma(L)$. 
Therefore, in this problem, it is sufficient to consider only links $L$ that satisfy $1-cr(L)\le\sigma(L)\le0$.
The $2$-bridge links show that the answer to this problem is yes for almost all pairs of integers. 
\begin{prop}\label{prop:2bridgesignature}
    For every pair of integers $(c,d)$ with $3-c\le d\le 0$ except for $(c,d)=(3,0),(5,0)$, there is a $2$-bridge link $L$ with $(cr(L),\sigma(L))=(c,d)$.
\end{prop}

In Section \ref{sec:examples}, we will provide specific examples of $2$-bridge links to support the proof of Proposition \ref{prop:2bridgesignature}.

Since the links $L$ satisfying $\sigma(L)=1-cr(L)$ are determined by Theorem \ref{thm:T(2,c)}, 
the next step is the classification of links satisfying $\sigma(L)=2-cr(L)$.
The main result of this paper is to determine all links with this property.

\begin{figure}[htb]
    \centering
    \includegraphics[scale=.3]{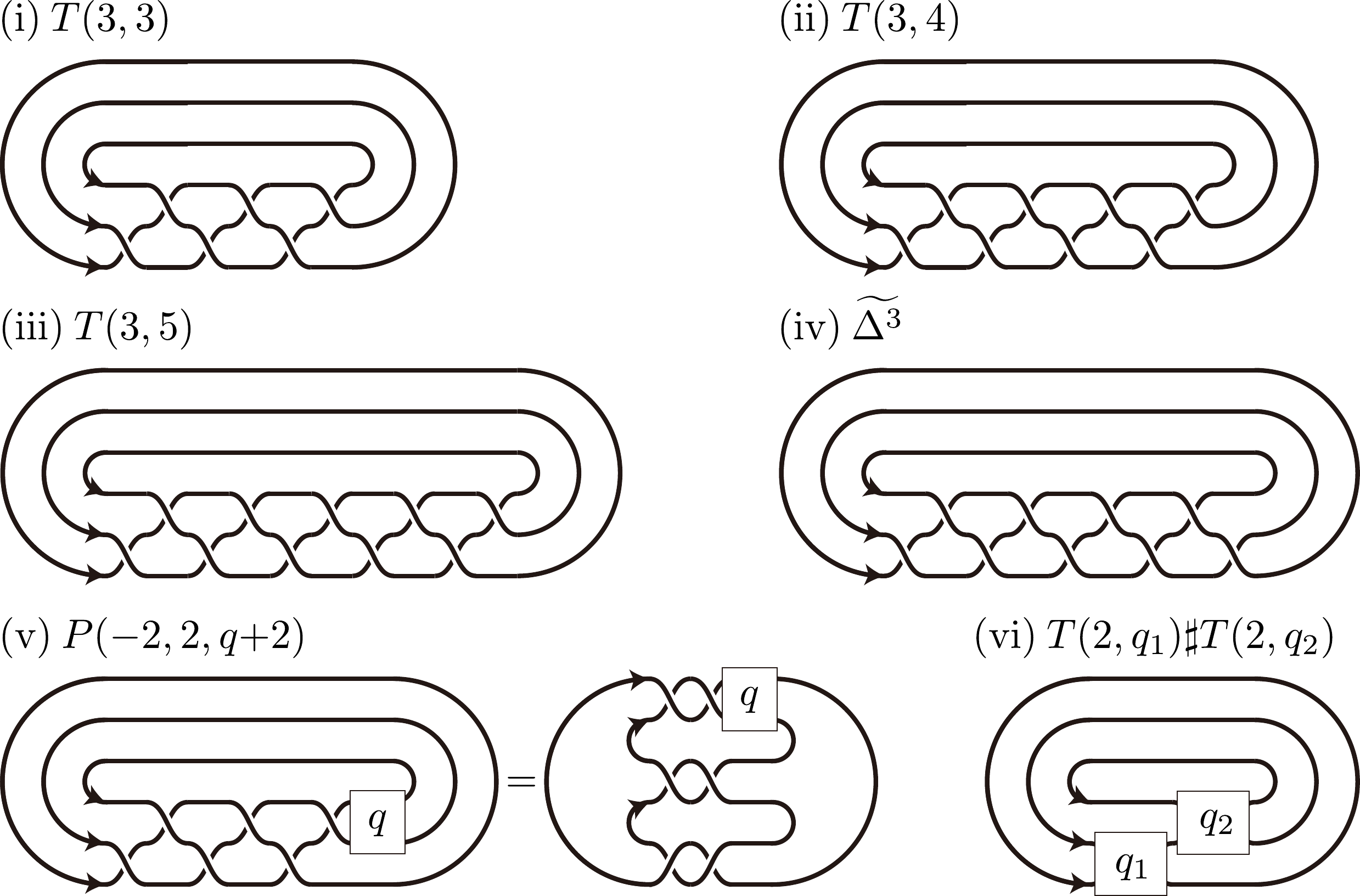}
    \caption{Oriented links $L$ with $\sigma(L)=2-cr(L)$.}
    \label{fig:Thm}
\end{figure}

\begin{thm}\label{thm:main}
Let $L$ be a non-split oriented link with $\sigma(L)=2-cr(L)\le 0$. 
Then $L$ is the closure of a positive $3$-braid and is one of the following.
\begin{enumerate}[\rm (i)]
\item $L=\widetilde{\Delta^2}=T(3,3)$. In this case, $cr(L)=6$ and $\sigma(L)=-4$. 
\item $L=\widetilde{\Delta^2\sigma_1\sigma_2}=T(3,4)$. In this case, $cr(L)=8$ and $\sigma(L)=-6$. 
\item $L=\widetilde{\Delta^2(\sigma_1\sigma_2)^2}=T(3,5)$. In this case, $cr(L)=10$ and $\sigma(L)=-8$. 
\item $L=\widetilde{\Delta^3}$. In this case, $cr(L)=9$ and $\sigma(L)=-7$.
\item $L=\widetilde{\Delta^2\sigma_2^q}=P(-2,2,q+2)\ (q>0)$. In this case, $cr(L)=q+6$ and $\sigma(L)=-q-4$.
\item $L=\widetilde{\sigma_1^{q_1}\sigma_2^{q_2}}=T(2,q_1)\sharp T(2,q_2)\ (q_1,q_2>1)$. In this case, $cr(L)=q_1+q_2$ and $\sigma(L)=-q_1-q_2+2$
\end{enumerate}

Here $\Delta$ is the $3$-braid $\sigma_1\sigma_2\sigma_1=\sigma_2\sigma_1\sigma_2$, $\widetilde{A}$ is the closure of a braid $A$, and $P(-2,2,p)$ is the oriented pretzel link of type $(-2,2,p)$, see Fig. \ref{fig:Thm}.
\end{thm}

\begin{rem}\label{rem:fiber}
All links in the conclusion of Theorem \ref{thm:main} are fibered, since they are
the closure of positive braids,  {\rm cf.}~\cite{Stallings}.
$T(3,4)$ and $T(3,5)$ are the only prime knots satisfying $\sigma(L)=2-cr(L)\le 0$.
\end{rem}

Combining Theorems \ref{thm:T(2,c)}, \ref{thm:main} and Proposition \ref{prop:2bridgesignature}, we have the following corollary that is the answer to Problem \ref{problem}.

\begin{figure}[htb]
    \centering
    \includegraphics[scale=.3]{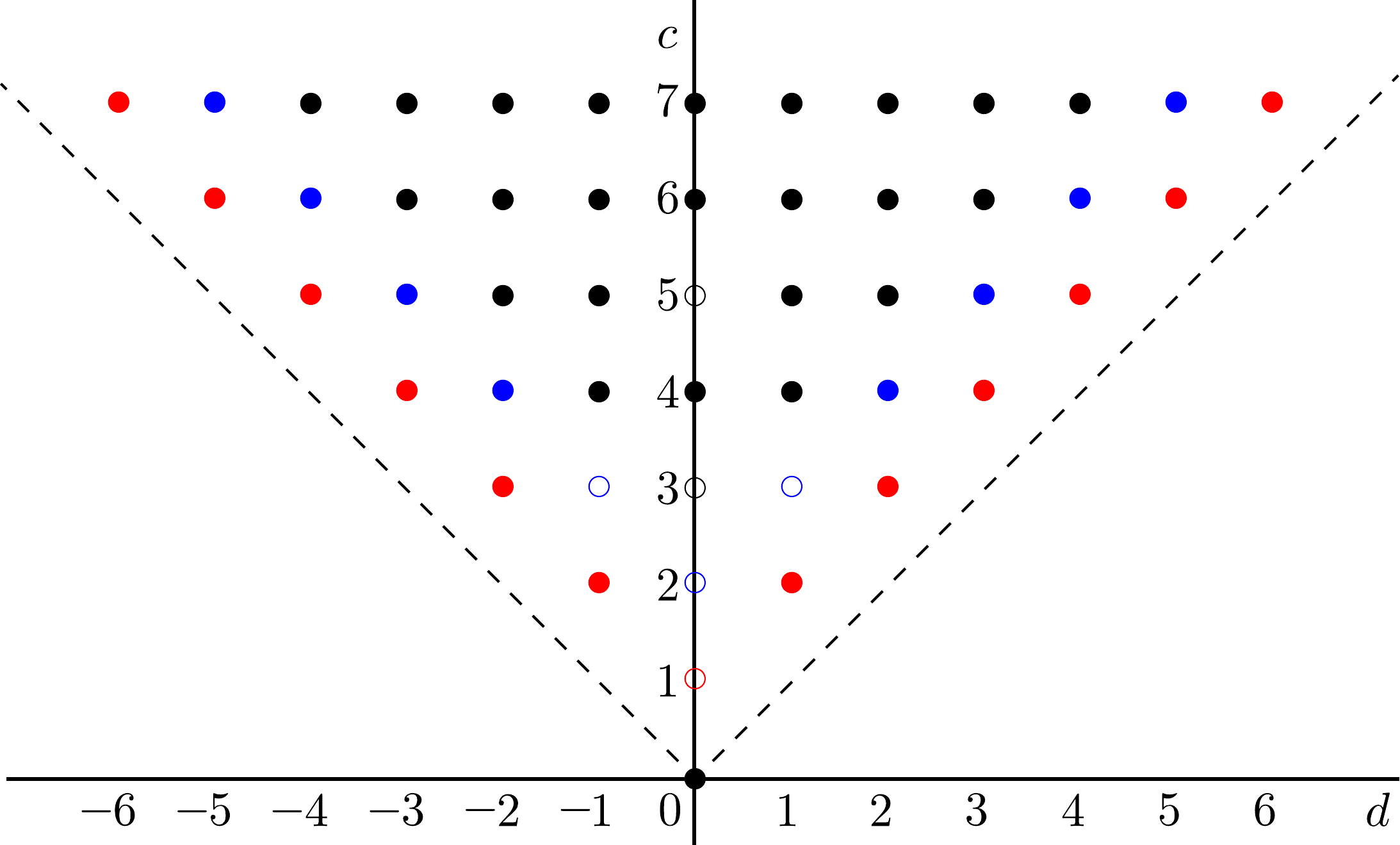}
    \caption{The colored circles $\bullet$ (resp. white circles $\circ$) represent the existence (resp. non-existence) of the links $L$ that satisfy $(cr(L),\sigma(L))=(c,d)$.
    The red circles correspond to $T(2,\pm c)$. The blue circles correspond to links in Theorem \ref{thm:main}.}
    \label{fig:table}
\end{figure}

\begin{cor}
     For every pair of integers $(c,d)$ with $1-c\le d\le c-1$ except for $(c,d)=(1,0),(2,0),(3,0),(3,\pm1),(5,0)$, there is an oriented link $L$ with $(cr(L),\sigma(L))=(c,d)$, see Fig. \ref{fig:table}.
\end{cor}

For the split sum $L_1\sqcup L_2$ of two links $L_1,L_2$, 
 $cr(L_1\sqcup L_2)=cr(L_1)+cr(L_2)$ and $\sigma(L_1)+\sigma(L_2)$ hold. 
 Obviously the trivial link $U_n$ with $n$-components satisfies $cr(U_n)=\sigma(U_n)=0$. 
Thus we have the following from Theorems \ref{thm:T(2,c)} and \ref{thm:main}:

\begin{cor}
Let $L$ be an oriented link. 
Then,  
$$-cr(L)\le\sigma(L)\le cr(L).$$
Moreover, the following hold.
\begin{enumerate}[\rm (1)]
    \item $\sigma(L)=cr(L)$ $(\sigma(L)=-cr(L))$ if and only if $L=U_n$.
    \item $\sigma(L)=1-cr(L)$ if and only if $L=T(2,c)$ or $L=T(2,c)\sqcup U_n$, where $c=cr(L)$.  
    \item $\sigma(L)=2-cr(L)$ if and only if $L=K$ or $L=K\sqcup U_n$ or $L=T(2,c_1)\sqcup T(2,c_2)$ or $L=T(2,c_1)\sqcup T(2,c_2)\sqcup U_n$, where $K$ is a link in Theorem \ref{thm:main} and $c_1,c_2$ are integers with $c_1,c_2>1,c_1+c_2=cr(L)$. 
\end{enumerate}
\end{cor}

One of the motivations of this study is a characterization of links obtained from $T(2,c)$ by a coherent band surgery.
We will discuss this application of Theorem \ref{thm:main} in Section \ref{sec:band}.


\section{Signature, crossing number, nullity and braid index}\label{sec:pre}
In this section, we refine the relation between the signature and crossing number in conjunction with nullity and braid index.
The {\em nullity} $n(L)$ is defined as that of $M+M^T$, where $M$ is a Seifert matrix of $L$ and $M^T$ is the transposed matrix of $M$. 
Note that $n(L)$ is bounded above by $\mu(L)-1$, where $\mu(L)$ is the number of components of $L$.
For a link diagram $D$ of $L$, let $s(D)$ be the number of Seifert circles of $D$. 
The braid index $br(L)$ is defined by the minimal number of strings among all braid diagrams of $L$.
Yamada showed that $br(L)= \min \{s(D)\mid D$ is a diagram of $L\}$ in \cite{Yamada}.

By examining the proof of Theorem \ref{thm:T(2,c)}, we can improve it as follows.

\begin{thm}\label{thm:ineq}
    Let $L$ be a non-trivial oriented knot or non-split oriented link with a diagram $D$. Then
\[
|\sigma(L)|+n(L)+s(D)\le cr(D)+1.
\]
    In particular,
\[
|\sigma(L)|+n(L)+br(L)\le cr(L)+1.
\]
\end{thm}

\begin{proof}
We assume that $L$ is a non-trivial oriented knot or non-split oriented link.
The canonical Seifert surface is obtained by attaching half-twisted bands to disks bounded by Seifert circles along all crossings. 
Then the canonical Seifert surface $F$, which is obtained from $D$, satisfies that 
$$\chi(F)=s(D)-cr(L),$$
where $\chi(F)$ is the Euler characteristic of $F$.
The absolute value of the signature of a matrix is at most the size of the matrix, and the size of a Seifert matrix is the first Betti number $b_1(F)$. 
Since $F$ is connected, $b_1(F)=1-\chi(F)=cr(D)-s(D)+1$. 
Hence, by this argument on the canonical Seifert surface $F$, we have the first inequality. 
$$|\sigma(L)|+n(L)\le cr(D)-s(D)+1$$
If we take $D$ as a minimal crossing diagram of $L$ so that $cr(D)=cr(L)$, we have the second inequality using the inequality $s(D)\ge br(L)$, that is, a result in \cite{Yamada}.
$$|\sigma(L)|+n(L)\le cr(L)-s(D)+1\le cr(L)-br(L)+1$$
This completes the proof of Theorem \ref{thm:ineq}.
\end{proof}

Theorem \ref{thm:T(2,c)} can be obtained from Theorem \ref{thm:ineq} as follows.
Since $L$ is not a trivial knot, we have $s>1$, thus
$$-cr(L)<\sigma(L)<cr(L).$$ 
Moreover, $s=2$ implies that the minimal diagram $D$ is a $2$-braid diagram of $L=T(2,c)$. 
Hence, $\sigma(L)=1-cr(L)$ if and only if $L=T(2,c)$, where $c=cr(L)$. 
\medskip

As a first step toward proving the main theorem (Theorem \ref{thm:main}),
we show that the link $L$ satisfying $\sigma(L)=2-cr(L)\le0$ is the closure of a positive $3$-braid with $cr(L)$ crossings.

\begin{prop}\label{prop:positive}
If a non-trivial oriented knot or non-split oriented link $L$ satisfies that $\sigma(L)=2-cr(L)\le0$, 
then $L$ is the closure of a positive $3$-braid with $cr(L)$ crossings. 
\end{prop}
To prove Proposition \ref{prop:positive} we use the following lemma.

\begin{lem}[\cite{M1}]\label{lem:murasugi}
Suppose an oriented link $L'$ is obtained from $L$ by smoothing on $c$ crossings. 
Then $\sigma(L')-c\le \sigma(L)\le\sigma(L')+c$. 
\end{lem}

\begin{figure}[htb]
    \centering
    \includegraphics[width=0.3\textwidth]{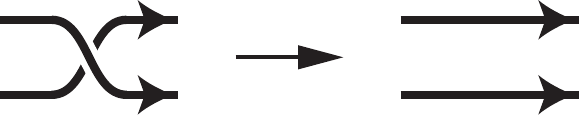}
    \caption{A smoothing on a crossing.}
    \label{fig:smoothing}
\end{figure}

\begin{proof}[Proof of Proposition \ref{prop:positive}]
By the same argument above, we have $s=3$, where $s$ is the number of Seifert circles on a minimal crossing diagram $D$ for $L$.
There are two cases for $D$: three Seifert circles nest or not. 
In the former case, $D$ is exactly a closed $3$-braid. 
In the latter case, $L$ is a connected sum of two closed $2$-braids, which is also a closed $3$-braid with $cr(L)$ crossings.

In both cases, we can assume that $L$ is the closure of a $3$-braid 
$A=\sigma_1^{p_1}\sigma_2^{q_1}\cdots\sigma_1^{p_r}\sigma_2^{q_r}$ with $cr(L)$ crossings. 
Let $p=\sum_{i=1}^{r}p_i$, $q=\sum_{i=1}^{r}q_i$,  $\overline{p}=\sum_{i=1}^{r}|p_i|$ and $\overline{q}=\sum_{i=1}^{r}|q_i|$. 
Note that $\overline{p},\overline{q}>1$ and $cr(L)=\overline{p}+\overline{q}$. 
By smoothing on the $\overline{q}-1$ crossings of $\sigma_2^{q_i}$'s,  we have a $(2,p)$-torus link $L'=T(2,p)$ from the link $L$.

There are three cases for $p$ as below:
If $p>0$, $\sigma(L')=1-p\ge1-\overline{p}$, and the equality holds only if $p_i>0$ for each $i$. 
If $p=0$, $\sigma(L')=0>1-\overline{p}$.
If $p<0$, $\sigma(L')=-1-p\ge0>1-\overline{p}$.
Then by Lemma \ref{lem:murasugi}, we have
$$2-cr(L)=2-\overline{p}-\overline{q}\le\sigma(L')-(\overline{q}-1)\le\sigma(L).$$ 

Since now $\sigma(L)=2-cr(L)$, $p_i>0$ for each $i$. 
By smoothing $\overline{p}-1$ crossings of $\sigma_1^{p_i}$'s, 
we also have $q_i>0$ for each $i$. 
This implies that the $3$-braid $A$ is positive. 
\end{proof}

In the next section, we will investigate the crossing number $cr(L)$ and the signature $\sigma(L)$, for the closure $L$ of a positive $3$-braid $A$. 
By Proposition \ref{prop:homogeneous} (\cite[Proposition 7.4]{M3}) below, in fact, $cr(L)$ generally agree with the number of crossings on $A$, except in the case where the braid index $br(L)$ is less than $3$.
A positive $n$-braid is a {\em homogeneous} $n$-braid, that is, an $n$-braid represented as $s_{i_1}s_{i_2}\cdots s_{i_c}$ with $1\le i_1,\ldots,i_c\le n-1$ and with $s_i=\sigma_i^{\pm1}$ for each $1\le i\le n-1$. Such a representation $s_{i_1}s_{i_2}\cdots s_{i_c}$ of homogeneous $n$-braid has $c$ crossings.
\begin{prop}[\cite{M3}]\label{prop:homogeneous}
    Suppose a link $L$ with braid index $n$ is the closure of a homogeneous $n$-braid with $c$ crossings, then $cr(L)=c$.
\end{prop}

\section{Positive $3$-braids}

The closed $3$-braids are well studied by Murasugi \cite{M2}. 
He classified all 3-braids into seven classes up to conjugation.
Recall that $\Delta$ is the $3$-braid $\sigma_1\sigma_2\sigma_1=\sigma_2\sigma_1\sigma_2$.
\begin{prop}[\cite{M2}]\label{prop:3-braid}
Any $3$-braid is conjugate to an element of the following sets $\Omega_0,\ldots, \Omega_6$.
\begin{enumerate}[\rm (1)]
\setcounter{enumi}{-1}
\item $\Omega_0=\{\Delta^{2n}\ |\ n\in \mathbb{Z}\}$.
\item $\Omega_1=\{\Delta^{2n}\sigma_1\sigma_2\ |\ n\in \mathbb{Z}\}$.
\item $\Omega_2=\{\Delta^{2n}(\sigma_1\sigma_2)^2\ |\ n\in \mathbb{Z}\}$.
\item $\Omega_3=\{\Delta^{2n+1}\ |\ n\in \mathbb{Z}\}$.
\item $\Omega_4=\{\Delta^{2n}\sigma_1^{-p}\ |\ n\in \mathbb{Z},p\in\mathbb{N}\}$.
\item $\Omega_5=\{\Delta^{2n}\sigma_2^q\ |\ n\in \mathbb{Z},q\in\mathbb{N}\}$.
\item $\Omega_6=\{\Delta^{2n}\sigma_1^{-p_1}\sigma_2^{q_1}\cdots\sigma_1^{-p_r}\sigma_2^{q_r}\ |\ n\in \mathbb{Z},p_i,q_i\in\mathbb{N}\}$.
\end{enumerate}
\end{prop}

Based on this classification, Murasugi further calculated the signatures of closed $3$-braids, which Erle later completed. 

\begin{prop}[\cite{M2,E}]\label{prop:signature3-braid}
Let $n$ be an integer, and let $p,q,p_i,q_i$ be positive integers. 
\begin{enumerate}[\rm (1)]
\setcounter{enumi}{-1}
\item For $A=\Delta^{2n}\in \Omega_0$, 
$\sigma(\widetilde{A})=-4n$.
\item For $A=\Delta^{2n}\sigma_1\sigma_2\in\Omega_1$, 
$\sigma(\widetilde{A})=
\begin{cases}
-4n&n\mbox{ is even}\\
-4n-2&n\mbox{ is odd}
\end{cases}$.
\item For $A=\Delta^{2n}(\sigma_1\sigma_2)^2\in\Omega_2$, 
$\sigma(\widetilde{A})=
\begin{cases}
-4n-2&n\mbox{ is even}\\
-4n-4&n\mbox{ is odd}
\end{cases}$.
\item For $A=\Delta^{2n+1}\in\Omega_3$, 
$\sigma(\widetilde{A})=
\begin{cases}
-4n-1&n\mbox{ is even}\\
-4n-3&n\mbox{ is odd}
\end{cases}$.
\item For $A=\Delta^{2n}\sigma_1^{-p}\in\Omega_4$, 
$\sigma(\widetilde{A})=
\begin{cases}
p-4n-1&n\mbox{ is even}\\
p-4n&n\mbox{ is odd}
\end{cases}$.
\item For $A=\Delta^{2n}\sigma_2^q\in\Omega_5$, 
$\sigma(\widetilde{A})=
\begin{cases}
-q-4n+1&n\mbox{ is even}\\
-q-4n&n\mbox{ is odd}
\end{cases}$.
\item For $A=\Delta^{2n}\sigma_1^{-p_1}\sigma_2^{q_1}\cdots\sigma_1^{-p_r}\sigma_2^{q_r}\in\Omega_6$, 
$\sigma(\widetilde{A})=\sum_{i=1}^{r}(p_i-q_i)-4n$.
\end{enumerate}
\end{prop}

To prove the main theorem (Theorem \ref{thm:main}), it is sufficient to determine $3$-braids in each family which are conjugate to positive $3$-braids.

\begin{prop}\label{prop:positive3-braid}~
\begin{enumerate}[\rm (1)]
\setcounter{enumi}{-1}
\item $\Delta^{2n}\in\Omega_0$ is conjugate to a positive $3$-braid if and only if $n\ge 0$.
\item $\Delta^{2n}\sigma_1\sigma_2\in\Omega_1$ is conjugate to a positive $3$-braid if and only if $n\ge0$.
\item $\Delta^{2n}(\sigma_1\sigma_2)^2\in\Omega_2$ is conjugate to a positive $3$-braid if and only if $n\ge0$.
\item $\Delta^{2n+1}\in\Omega_3$  is conjugate to a positive $3$-braid if and only if $n\ge0$.
\item $\Delta^{2n}\sigma_1^{-p}\in\Omega_4$  is conjugate to a positive $3$-braid if and only if $n\ge \frac{p}{2}>0$
\item $\Delta^{2n}\sigma_2^q\in\Omega_5$ is conjugate to a positive $3$-braid if and only if $n\ge0$.
\item $\Delta^{2n}\sigma_1^{-p_1}\sigma_2^{q_1}\cdots\sigma_1^{-p_r}\sigma_2^{q_r}\in\Omega_6$ is conjugate to a positive $3$-braid if and only if $n\ge \frac{p}{2}>0$, 
where $p=\sum_{i=1}^{r}p_i$.
\end{enumerate}
\end{prop}

Before proving Proposition \ref{prop:positive3-braid}, let us make some observations about the crossing number to show the main theorem (Theorem \ref{thm:main}). 
In the cases of (0), (1), (2), (3), and (5) in Proposition \ref{prop:positive3-braid}, it is easily verified that the $3$-braids 
$\Delta^{2n}$, $\Delta^{2n}\sigma_1\sigma_2$, $\Delta^{2n}(\sigma_1\sigma_2)^2$, $\Delta^{2n+1}$, $\Delta^{2n}\sigma_2^q$ for $n\ge0$ 
are in fact positive $3$-braids and 
\begin{align*}
    cr(\Delta^{2n})&=6n,\\cr(\Delta^{2n}\sigma_1\sigma_2)&=6n+2,\\cr(\Delta^{2n}(\sigma_1\sigma_2)^2)&=6n+4,\\cr(\Delta^{2n+1})&=6n+3,\\cr(\Delta^{2n}\sigma_2^q)&=6n+q.
\end{align*}
Here, $cr(A)$ represents the number of crossings in a positive braid $A$. 
In the other cases (4) and (5) in Proposition \ref{prop:positive3-braid}, $\Delta^{2n}\sigma_1^{-p}$ and $\Delta^{2n}\sigma_1^{-p_1}\sigma_2^{q_1}\cdots\sigma_1^{-p_r}\sigma_2^{q_r}$ for $n\ge \frac{p}{2}>0$ are also positive $3$-braids themselves, and 
\begin{align*}
    cr(\Delta^{2n}\sigma_1^{-p})&=6n-p,\\
    cr(\Delta^{2n}\sigma_1^{-p_1}\sigma_2^{q_1}\cdots\sigma_1^{-p_r}\sigma_2^{q_r})&=6n+\sum_{i=1}^{r}(q_i-p_i),
\end{align*}
as $\Delta\sigma_i^{\pm1}=\sigma_j^{\pm1}\Delta$ and $\Delta\sigma_i^{-1}=\sigma_i\sigma_j$ for $\{i,j\}=\{1,2\}$. 
Recalling Proposition \ref{prop:homogeneous}, 
in most of cases $cr(L)=cr(A)$ for $L=\widetilde{A}$, but $cr(L)<cr(A)$ in special cases where $br(L)<3$.
Such exceptions in Proposition \ref{prop:positive3-braid} are, the case where $n=0$ for (1), (2), (3), $n=p=1$ for (4), and $n=r=p_1=1$ for (6). 
In these cases, 
$A=\sigma_1\sigma_2$, $(\sigma_1\sigma_2)^2=\sigma_1^2\sigma_2\sigma_1$, 
$\Delta=\sigma_1\sigma_2\sigma_1$, $\Delta^2\sigma_1^{-1}=\sigma_2^2\sigma_1\sigma_2^2$, $\Delta^2\sigma_1^{-1}\sigma_2^{q_1}=\sigma_2^2\sigma_1\sigma_2^{q_1+2}$, thus 
$L$ is the unknot $0_1$, the trefoil knot $T(2,3)$, the Hopf link $T(2,2)$, $T(2,4)$, $T(2,q_1+4)$, respectively.
Hence $br(L)<3$ and $\sigma(L)\le 1-cr(L)$.
Using this observation, we prove Theorem \ref{thm:main}. 
\begin{proof}[Proof of Theorem \ref{thm:main}]
    Suppose that a non-trivial oriented knot or non-split oriented link $L$ satisfies $\sigma(L)=2-cr(L)\le 0$. 
    By Propositions \ref{prop:positive}, \ref{prop:3-braid} and \ref{prop:positive3-braid}, $L$ is the closure of a positive $3$-braid $A\in \Omega_i$ for some $i\in \{0,1,\ldots,6\}$ that appear in Proposition \ref{prop:positive3-braid}.
    From Proposition \ref{prop:signature3-braid} and the above observation, it can be easily seen that the condition $\sigma(L)=2-cr(L)$ is satisfied only if $n=1$ in each case of (0), (2), (3) and (5) in Proposition \ref{prop:positive3-braid}. 
    In the case of (4), the condition $\sigma(L)=2-cr(L)$ is satisfied only if $n=1$ and $p=2$. 
    In this case, $A=\Delta^2\sigma_1^{-2}=\sigma_2\sigma_1^2\sigma_2$ and $L=T(2,2)\sharp T(2,2)$.     
    In the case of (6), similarly, the condition $\sigma(L)=2-cr(L)$ is satisfied only if $n=p_1=p_2=1$ and $r=2$. 
    In this case, $A=\Delta^2\sigma_1^{-1}\sigma_2^{q_1}\sigma_1^{-1}\sigma_2^{q_2}=\sigma_2\sigma_1^{q_1+2}\sigma_2^{q_2+1}$ and $L=T(2,q_1+2)\sharp T(2,q_2+2)$. 
    This completes the proof of Theorem \ref{thm:main}.
\end{proof}

To prove Proposition \ref{prop:positive3-braid}, we apply the algorithm introduced by Elrifai and Morton \cite{EM} to $3$-braids. 
Here we briefly review their argument.

Let $B_3$ be the set of all $3$-braids, and let $B_3^+$ be the set of all positive $3$-braids including the trivial braid $e$.
An automorphism $\tau:B_3\to B_3$ is defined by $\tau(\sigma_1)=\sigma_2$ and $\tau(\sigma_2)=\sigma_1$.
A partial order $\le$ on $B_3$ is defined as follows:
For $A,B\in B_3$, $A\le B$ if $B=C_1AC_2$ for some $C_1,C_2\in B_3^+$. 
By definition, $B\in B_3^+$  
if and only if $e\le B$. 
Since $C\Delta=\Delta\tau(C)$ for any $C\in B_3$, if $\Delta^r\le B\le\Delta^s$ then $B=E_1\Delta^r=E_2\Delta^r$ 
and $\Delta^s=D_1B=BD_2$ for some $D_1,D_2,E_1,E_2\in B_3^+$, see \cite[Propositions 1.2 and 1.3]{EM}. 
Moreover the following theorem (\cite[Theorem 1.5]{EM}) is shown.
\begin{thm}[\cite{EM}]
Every $B$ satisfies $\Delta^r\le B\le\Delta^s$ for some $r,s\in\mathbb{Z}$.
\end{thm}
The subset $\{B\in B_3\ |\ \Delta^r\le B\le\Delta^s\}$ of $B_3$ is denoted by the interval $[r,s]$.
For $B\in B_3$, $\inf B$ and $\sup B$ is defined as follows: 
$$\inf B=\max\{r \mid \Delta^r\le B\}, \quad
\sup B=\min\{s \mid B\le\Delta^s\}.$$
If $B\ge \Delta^r$, put $P=\Delta^{-r}B\ge e$.
Then we have $\inf(B)=\inf(P)+r$ and $\sup(B)=\sup(P)+r$.
For $P\ge e$, {\em starting set} $S(P)$ and the {\em finishing set} $F(P)$ are defined as follows:
$$S(P)=\{i\mid P=\sigma _i P_i,P_i\geq e\},\quad
F(P)=\{ i\mid P=P_i\sigma _i,P_i\geq e\}.
$$
It was shown in \cite[Theorem 2.9]{EM} that 
every positive braid $P$ has the unique expression $P = P_1P_2\cdots P_k$, which is called the {\em left-canonical form},
satisfying $P_i\in [0,1]$, $P_k\neq e$ and $S(P_{i+1})\subset F(P_i)$ for each $i$. 
From the left-canonical form, we can find $\inf P$ and $\sup P$ as follows,  see \cite[Theorem 2.11]{EM}:
$$\inf P=\max\{i\ |\ P_i=\Delta \}, \quad \sup P=k.$$
Let $A = \Delta^r P_1P_2\ldots P_k$, where $P_1\cdots P_k$ is the left-canonical form of the positive braid $\Delta^{-r}A$ and $P_1\neq \Delta$. 
Then $r=\inf A$ and $r+k=\sup A$. 
$c(A) = \Delta^r P_2\cdots P_k\tau^r(P_1)$ is said to be given by {\em cycling} $A$.
Note that $c(A)$ is conjugate to $A$ since $\tau^r(P_1)c(A)=A\tau^r(P_1)$.  
The following lemma allows us to verify whether or not a given braid $A$ can be conjugated to a braid $B$ with $\inf B>\inf A$, see \cite[Lemma 4.3]{EM}.

\begin{lem}[\cite{EM}]\label{lem:4.3}
Suppose that $A$ is conjugate to $B$ with $\inf B>\inf A$. 
Then repeated cycling will produce $c^j(A)$ with $\inf c^j(A ) > \inf A$, for some $j$.
\end{lem}
Applying \cite[Theorem 2.6 and Lemma 2.7]{EM} to $3$-braids we can verify the following. 

\begin{lem}\rm~
\begin{enumerate}
\item 
$[0,1]=\{e,\sigma_1,\sigma_2,\sigma_1\sigma_2,\sigma_2\sigma_1,\Delta\}$.
\item
$S(\sigma_1)=S(\sigma_1\sigma_2)=F(\sigma_1)=F(\sigma_2\sigma_1)=\{1\}$, \\
$S(\sigma_2)=S(\sigma_2\sigma_1)=F(\sigma_2)=F(\sigma_1\sigma_2)=\{2\}$, 
and $S(\Delta)=F(\Delta)=\{1,2\}$.
\end{enumerate}
\end{lem}

By using these arguments, we give the proof of Proposition \ref{prop:positive3-braid}.

\begin{proof}[Proof of Proposition \ref{prop:positive3-braid}]
For (0), (1), (2), (3) of  Proposition \ref{prop:positive3-braid},
we can see $\Delta^{2n}$, $\Delta^{2n}\sigma_1\sigma_2\in[2n,2n+1]$
and $\Delta^{2n}(\sigma_1\sigma_2)^2=\Delta^{2n+1}\sigma_2$, $\Delta^{2n+1}\in [2n+1,2n+2]$. 
Since two braids contained in disjoint intervals cannot be conjugate, 
these $3$-braids cannot be conjugate to a $3$-braid $B$ with $\Delta^{2n+2}\le B$.
\begin{enumerate}[\rm (1)]
\setcounter{enumi}{3}
\item Let $A=\Delta^{2n}\sigma_1^{-p}\in\Omega_4$, and put $P=\Delta^{p-2n}A$,
where $n$ is an integer and $p$ is a positive integer. 
\begin{align*}
    P=\Delta^{p}\sigma_1^{-p}&=\big(\Delta\tau^{p-1}(\sigma_1^{-1})\big)\big(\Delta\tau^{p-2}(\sigma_1^{-1})\big)\cdots\big(\Delta\sigma_1^{-1}\big)\\
    &=\big(\tau^{p-1}(\sigma_1\sigma_2)\big)\big(\tau^{p-2}(\sigma_1\sigma_2)\big)\cdots\big(\sigma_1\sigma_2\big)
\end{align*}
This is the left-canonical form of $P$, and then $\inf A=2n-p$ and $\sup A=2n$ since $A=\Delta^{2n-p}P$. 
By cycling $A$, we have 
\begin{align*}
c(A)&=\Delta^{2n-p}\big(\tau^{p-2}(\sigma_1\sigma_2)\big)\big(\tau^{p-3}(\sigma_1\sigma_2)\big)\cdots\big(\sigma_1\sigma_2\big)\big(\tau^{2n-1}(\sigma_1\sigma_2)\big)\\
&=\Delta^{2n-p}\big(\tau^{p-1}(\sigma_2\sigma_1)\big)\big(\tau^{p-2}(\sigma_2\sigma_1)\big)\cdots\big(\tau(\sigma_2\sigma_1)\big)\big(\sigma_2\sigma_1\big),
\end{align*}
and still $\inf c(A)=2n-p$ and $\sup c(A)=2n$.
Since $c^2(A)=A$, by Lemma \ref{lem:4.3}, $A$ cannot be conjugate to a $3$-braid $B$ with $\Delta^{2n-p+1}\le B$. 
\item Let $A=\Delta^{2n}\sigma_2^q\in\Omega_5$, where $n$ is an integer and $q$ is a positive integer. 
Since $(\sigma_2)^q$ is the left-canonical form for $\Delta^{-2n}A$, 
$\inf A=2n$ and $\sup A=2n+q$.  
Since $c(A)=A$, by Lemma \ref{lem:4.3}, $A$ cannot be conjugate to a $3$-braid $B$ with $\Delta^{2n+1}\le B$.
\item Let $A=\Delta^{2n}\sigma_1^{-p_1}\sigma_2^{q_1}\cdots\sigma_1^{-p_r}\sigma_2^{q_r}\in\Omega_6$, where $p_i,q_i$ are positive integers. 
By replacing $\sigma_1^{-p_i}$ with $\sigma_1^{-1}\sigma_2^0\sigma_1^{-1}\cdots\sigma_2^0\sigma_1^{-1}$ if $p_i>1$, 
we can rewrite $A$ as \\ $\Delta^{2n}\sigma_1^{-1}\sigma_2^{q'_1}\sigma_1^{-1}\cdots\sigma_1^{-1}\sigma_2^{q'_p}$, where $p=\sum_{i=1}^{r}p_i$ and $q_i'\ge0$ for each $i$, especially $q'_p=q_r>0$.
 Put $P=\Delta^{p-2n}A$, then we have 
\begin{align*}
P&=
\big(\Delta\tau^{p-1}(\sigma_1^{-1}\sigma_2^{q'_1})\big)\big(\Delta\tau^{p-2}(\sigma_1^{-1}\sigma_2^{q'_2})\big)\cdots
\big(\Delta\sigma_1^{-1}\sigma_2^{q'_p}\big)\\
&=\big(\tau^{p-1}(\sigma_1\sigma_2)\big)\big(\tau^{p-1}(\sigma_2)\big)^{q'_1}\big(\tau^{p-2}(\sigma_1\sigma_2)\big)\big(\tau^{p-2}(\sigma_2)\big)^{q'_2}\cdots
\\&\hskip 90pt \cdots
\big(\tau(\sigma_1\sigma_2)\big)\big(\tau(\sigma_2)\big)^{q'_{p-1}}
\big(\sigma_1\sigma_2\big)\big(\sigma_2\big)^{q'_p}.
\end{align*}
This is the left-canonical form of $P$, and then $\inf A=2n-p$ and $\sup A=2n+q$ since $A=\Delta^{2n-p}P$, 
where $q=\sum_{i=1}^{p}q'_i=\sum_{j=1}^{r}q_j$. 
By the cycling $A$, we have 
\begin{align*}
  c(A)=
\Delta^{2n-p}\big(\tau^{p-1}(\sigma_2)\big)^{q'_1}\big(\tau^{p-2}(\sigma_1\sigma_2)\big)\big(\tau^{p-2}(\sigma_2)\big)^{q'_2}\cdots
\\\hspace{100pt}\cdots\big(\sigma_1\sigma_2\big)\big(\sigma_2\big)^{q'_p}\big(\sigma_2\sigma_1\big).  
\end{align*}
For $1\le j\le q'_1$, by additional $j$ time cyclings, we have
\begin{align*}
    c^{j+1}(A)&=
\Delta^{2n-p}\big(\tau^{p-1}(\sigma_2)\big)^{q'_1-j}\big(\tau^{p-2}(\sigma_1\sigma_2)\big)\big(\tau^{p-2}(\sigma_2)\big)^{q'_2}\cdots
\\&\hspace{100pt}\cdots
\big(\sigma_1\sigma_2\big)\big(\sigma_2\big)^{q'_p}\big(\sigma_2\sigma_1\big)\big(\sigma_1\big)^{j},\\
c^{q'_1+1}(A)&=
\Delta^{2n-p}\big(\tau^{p-2}(\sigma_1\sigma_2)\big)\big(\tau^{p-2}(\sigma_2)\big)^{q'_2}\cdots
\\&\hspace{100pt}\cdots
\big(\sigma_1\sigma_2\big)\big(\sigma_2\big)^{q'_p}\big(\sigma_2\sigma_1\big)\big(\sigma_1\big)^{q'_1}.
\end{align*}
Finally, we have
\begin{align*}
c^{p+q}(A)&=
\Delta^{2n-p}\big(\sigma_2\sigma_1\big)\big(\sigma_1\big)^{q'_1}\big(\tau(\sigma_2\sigma_1)\big)\big(\tau(\sigma_1)\big)^{q'_2}\cdots
\\&\hspace{100pt}\cdots
\big(\tau^{p-1}(\sigma_2\sigma_1)\big)\big(\tau^{p-1}(\sigma_1)\big)^{q'_p}\\
&=\tau^{p}(A).
\end{align*}
Since $c^j(\tau^p(A))=\tau^j(c^j(A))$, $c^{2p+2q}(A)=c^{p+q}(\tau^p(A))=A$, by Lemma \ref{lem:4.3}, $A$ cannot be conjugate to a $3$-braid $B$ with $\Delta^{2n-p+1}\le B$.
\end{enumerate}
\end{proof}


\section{Signature and crossing number of $2$-bridge knots and links}\label{sec:examples}
In this section, we provide specific examples of $2$-bridge links to support the proof of Proposition \ref{prop:2bridgesignature}. 
All $2$-bridge links are alternating, that is, having alternating diagrams. 
It is well known that the crossing number of a reduced alternating diagram of an alternating link is the crossing number of the link, {\rm cf}. \cite{Kauffman, TaitCoj, Thistlethwaite}. 
On the other hand, there are several methods for calculating the signature of $2$-bridge knot, {\rm cf.} \cite{GJ, M4}. 

\begin{figure}[htb]
    \centering
    \includegraphics[scale=.3]{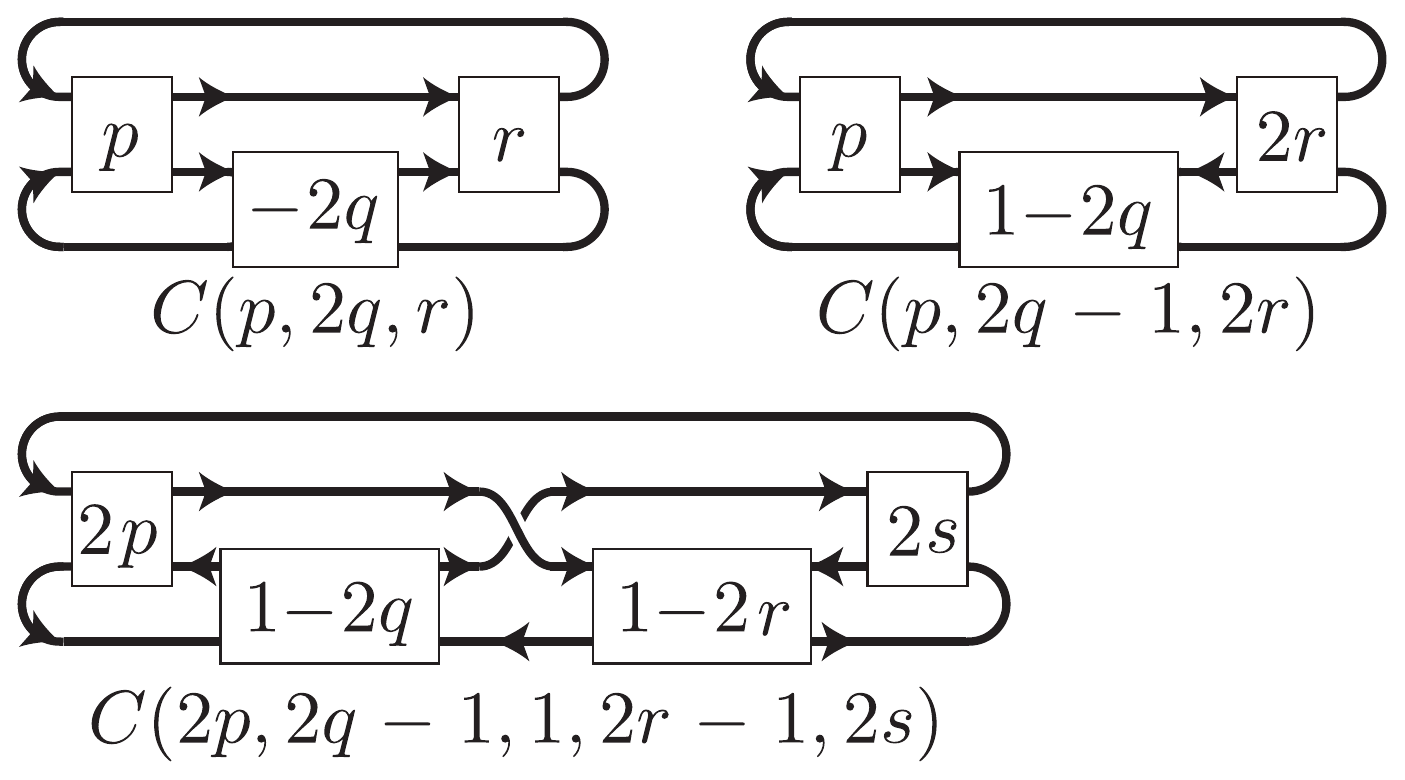}
    \caption{Oriented $2$-bridge links  
    $C(p,2q,r)$, $C(p,2q-1,2r)$, $C(2p,2q-1,1,2r-1,2s)$.}
    \label{fig:examples}
\end{figure}

We consider $2$-bridge links, which are denoted by Conway notations $C(p,2q,r)$, $C(p,2q-1,2r)$ and $C(2p,2q-1,1,2r-1,2s)$ where $p,q,r,s$ are positive integers and the diagrams are oriented and alternating as depicted in Fig. \ref{fig:examples}.
Note that $C(2p,2q-1,1,2r-1,2s)$ is a knot for each $p, q,r,s$, 
and $C(p,2q,r)$ (resp. $C(p,2q-1,2r)$) is a knot if $p+r$ (resp. $p$) is odd. Otherwise, it is a $2$-component link.
For these $2$-bridge links, we have the following.
\begin{lem}\label{lem:2bridge} 
\begin{align}
    cr(C(p,2q,r))&=p+2q+r,\\
    \sigma(C(p,2q,r))&=1-p-r\notag\\
    cr(C(p,2q-1,2r))&=p+2(q+r)-1,\\
    \sigma(C(p,2q-1,2r))&=1-p\notag\\
    cr(C(2p,2q-1,1,2r-1,2s))&=2(p+q+r+s)-1,\\
    \sigma(C(2p,2q-1,1,2r-1,2s))&=0\notag
\end{align}
\end{lem}

\begin{proof} 
    Since Fig. \ref{fig:examples} shows reduced and alternating diagrams, the crossing number is easily seen for each. 
    So let us consider the signature for each of the links.
    In the case of knots, we can easily verify the conclusion using a method in \cite{GJ}.
    In the case of links, we only consider the link $L=C(p,2q-1,2r)$ where $p$ is even, but the same argument can be made for the link $C(p,2q,r)$ where $p+r$ is even.     
    $L$ is obtained from $L_1=C(p+1,2q-1,2r)$ by a smoothing on a crossing, and similarly $L_2=C(p-1,2q-1,2r)$ is obtained from $L$ by a smoothing on a crossing. 
    $L_1$ and $L_2$ are knots with $\sigma(L_1)=-p$, $\sigma(L_2)=2-p$, thus, 
    by Lemma \ref{lem:murasugi}, we have $\sigma(L)=1-p$. 
    $$1-p=\sigma(L_2)-1\le \sigma(L)\le \sigma(L_1)+1=1-p$$  
\end{proof}
Then we can give a proof of Proposition \ref{prop:2bridgesignature} using Lemma \ref{lem:2bridge}.
\begin{proof}[Proof of Proposition \ref{prop:2bridgesignature}]
    Let $(c,d)$ be a pair of integers with $3-c\le d\le 0$, 
    other than $(3,0)$ and $(5,0)$. 
    In the case where $c+d$ is odd and $d<0$, we can choose positive integers $p,q,r$ so that $p+2q+r=c$, $1-p-r=d$. 
    Then the link $L=C(p,2q,r)$ satisfies $cr(L)=c$ and $\sigma(L)=d$. 
    In the case where $c$ is odd and $d=0$, since $c\ge7$, we can chose positive integers $p,q,r,s$ so that $2(p+q+r+s)-1=c$. 
    Then the link $L=C(2p,2q-1,1,2r-1,2s)$ satisfies $cr(L)=c$ and $\sigma(L)=0$. 
    In the case where $c+d$ is even, we can chose positive integers $p,q,r$ so that $p+2(q+r)-1=c$, $1-p=d$. 
    Then the link $L=C(p,2q-1,2r)$ satisfies $cr(L)=c$ and $\sigma(L)=d$. 
\end{proof}

\section{Applications}\label{sec:band}

\begin{figure}[htb]
    \centering
    \includegraphics[width=0.5\textwidth]{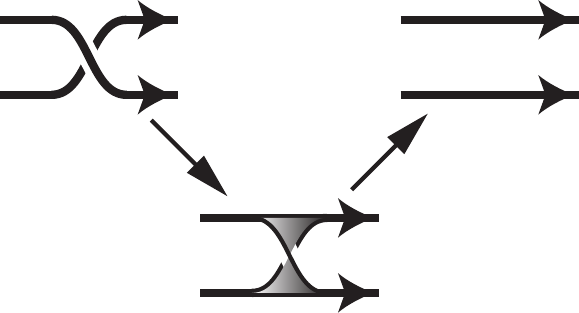}
    \caption{A smoothing can be considered as a coherent band surgery.}
    \label{fig:band}
\end{figure}

In this section, we discuss applications of the results to band surgery problems.
For an oriented link $L$ in $S^3$,
suppose $b: I\times I \to S^3$ is an embedding such that $b(I\times I)\cap L=b(I \times \partial I)$, where $I$ is a closed interval $[0,1]$.
Then we say that $L'=(L-b(I\times \partial I))\cup b(\partial I \times I)$ is obtained from $L$ by a band surgery along $b$.
A band surgery is said to be {\em coherent} if the orientations of $L$ and $L'$ coincide outside $b$.
A smoothing and a coherent band surgery can be considered as the same operation, see Fig. \ref{fig:band}.

In \cite{NatPhys}, it is shown that a vortex knot of trefoil knot $T(2,3)$ 
is changed into one of Hopf link $T(2,2)$ by an anti-parallel reconnection.
In \cite{vortex}, smoothing is used to model an anti-parallel reconnection 
and shows that the helicity is preserved under an anti-parallel reconnection.

In \cite{PNAS, SciRep}, it is shown that the XerCD-{\it dif}-FtsK system unlinks replication DNA catenanes in a stepwise manner.
The link types of the replication catenanes are torus link $T(2,c)$.
If we assign orientations for each component by using a sequence of the {\it dif} site, they have parallel orientation and $\sigma(T(2,c))=1-c$.
The recombination is modeled by using a smoothing (or a coherent band surgery).
In \cite{PNAS}, by using Theorem \ref{thm:T(2,c)}, it is shown that a knot or a link $L'$ is obtained from $L=T(2,c)$ by a smoothing and $cr(L')<cr(L)$, then $L'=T(2,c-1)$.
As an application, we can uniquely characterize the unlinking pathway of the replication catenane of type $T(2,c)$ if we assume that the crossing number goes down at each recombination.
In \cite{SciRep}, we relaxed the assumption a little and considered the case where $cr(L')=cr(L)$ for a smoothing on $L=T(2,c)$. 
Then $L'$ satisfies $\sigma(L')=2-cr(L')$.
This is one of the motivations for characterizing knots and links with this property.

\begin{figure}[htb]
    \centering
    \includegraphics[width=0.5\textwidth]{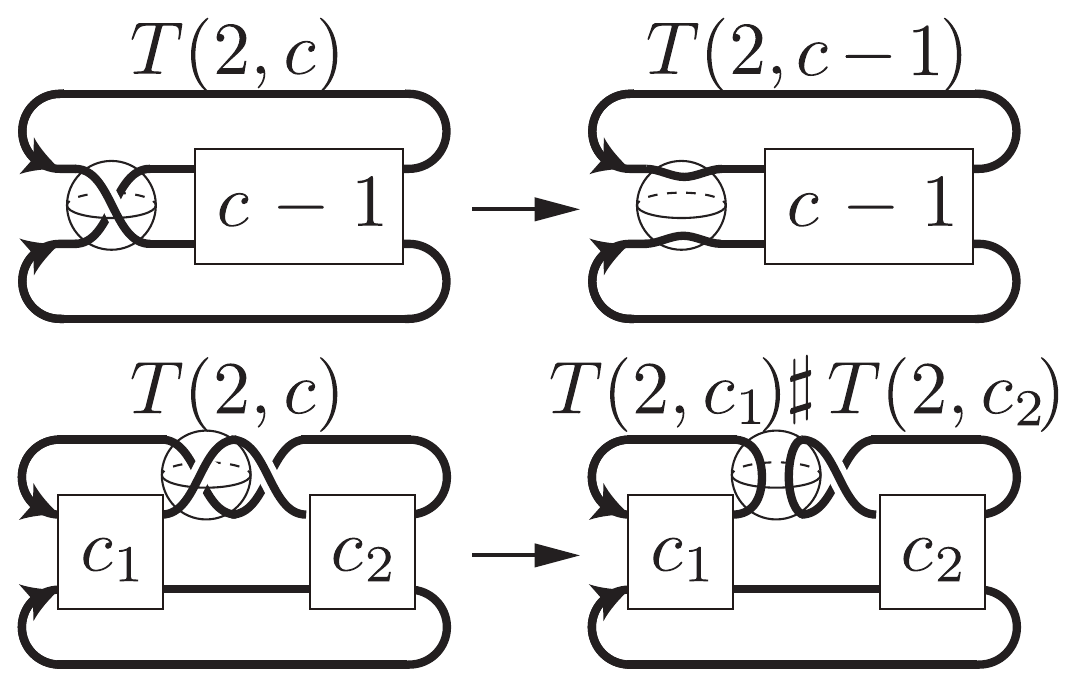}
    \caption{Smoothings on $T(2,c)$.}
    \label{fig:T(2c)smoothing}
\end{figure}

On the other hand, an equivalence relation can be introduced to characterize the smoothings.
Let us consider two smoothings from $L_i$ to $L_i'$ ($i=1,2$) in $S^3$, and let $B_i$ be a $3$-ball in $S^3$ that contains the crossing of $L_i$ as a tangle where the smoothing occurs. Thus, it is considered that $L_i$ and $L_i'$ coincide outside of $B_i$ for each $i\in\{1,2\}$.
We say these smoothings are equivalent if there exists an orientation preserving homeomorphism $h:S^3\to S^3$ such that $h(L_1)=L_2$, $h(L_1')=L_2'$ and $h(B_1)=B_2$.
Smoothings between fibered knots and links are characterized with respect to this equivalence relation in \cite{JLMS}.
If a fibered knot or link $L'$ is obtained from $L=T(2,c)$ by a smoothing and $\chi(L')>\chi(L)$, then $L'$ is either (1) $L'=T(2,c-1)$, or (2) $L'=T(2,c_1)\# T(2,c_2)$ where $c_1,c_2>1, c_1+c_2=c$. Furthermore, it was shown that there is a unique smoothing up to equivalence for each link in (1) and (2) determined as $L'$. 
Here $\chi(L)$ denotes the Euler characteristic of $L$, that is, the maximal value of the Euler characteristic of all (possibly disconnected) Seifert surfaces for $L$.
As stated in Remark \ref{rem:fiber}, the links in the conclusion of Theorem \ref{thm:main} are fibered. 
According to \cite{JLMS}, we can demonstrate that the above conclusion holds even without assuming that $L'$ is fibered, cf. Theorem \ref{thm:T(2c)smoothing}.
This fact has already been described in Theorem 21 and Corollary 22 in 
\cite{Ishihara}.
\begin{thm}\label{thm:T(2c)smoothing} Let $L=T(2,c)$ and $L'$ an oriented link with $cr(L')\le c$ for $c>1$.
    Suppose $L'$ is obtained from $L$ by smoothing a crossing. Then $L'=T(2,c-1)$ or $T(2,c_1)\# T(2,c_2)$ where $c_1,c_2>1, c_1+c_2=c$. Moreover, the smoothing is unique up to equivalence for each $L'$, see Fig. \ref{fig:T(2c)smoothing}.
\end{thm}
\begin{rem}
    For any Seifert surface $F$ of $L$, we have an inequality $|\sigma(L)|\le b_1(F)$, where $b_1(F)$ is the first Betti number of $F$, because the size of Seifert matrix obtained from $F$ is $b_1(F)$. Suppose $F$ consists of $\mu$ components, then $b_1(F)=\mu-\chi(F)$, not $1-\chi(F)$.  
    In the case where $L$ has a disconnected taut Seifert surface $F$, which maximises the Euler characteristic ($\chi(F)=\chi(L)$), the inequality $|\sigma(L)|\le 1-\chi(L)$ in Proposition 8 in \cite{Ishihara} may not be correct. 
    For a fibered link $L$, however, a fiber surface is known as a unique taut Seifert surface that is connected, thus the inequality is correct.
\end{rem}

\section*{Acknowledgment}
The authors would like to thank Kouki Taniyama for the discussion on the pair of knot invariants.
The first author is supported by JSPS Kakenhi grant number JP23K03114.
The last author is supported by JSPS Kakenhi grant number JP23K17652.
 Part of this research was conducted at the Banff International Research Station. We would like to thank them for their hospitality.


\end{document}